\documentclass[12pt]{article}

\usepackage{amsmath}
\usepackage{amssymb}
\usepackage{amsthm}
\usepackage{graphicx}
\usepackage{booktabs}
\usepackage[ruled,vlined,linesnumbered]{algorithm2e}
\usepackage{enumitem}
\usepackage{hyperref}

\newtheorem{theorem}{Theorem}

\newtheorem{propo}{Proposition}

\title{Acceleration of the Method of Alternating Projections on subspaces by GMRES}

\title{Acceleration of the Method of Alternating Projections on subspaces by GMRES}

\author{Luis Manuel Hernández-Ramos,\and Maryam Calles-Cabriles,\and Gorka Emaldi-Porco.
  \thanks{%
    Escuela de Computación. Facultad de Ciencias.\\
    Universidad Central de Venezuela.\\
    \texttt{email: luis.m.hernandez@ucv.ve}}}

\date{\today}

\date{\today}

\newcommand{\keywords}[1]{%
  \par\noindent\textbf{Keywords:} #1\par\medskip
}

\begin{document}

\maketitle
\thispagestyle{empty}

\begin{abstract}
In this paper we propose a GMRES-accelerated version of the Method of Alternating 
Projections (MAP) for subspaces. It is proved that the accelerated version not only 
converges to a feasible solution, but also solves a best approximation problem. 
Comparisons are made with other accelerations proposed in the literature, for solving 
saddle point problems.
\end{abstract}


\section*{Introduction}
The Method of Alternating Projections (MAP) is an algorithm for finding the best 
approximation to any given point in a Hilbert space from the intersection of a 
finite number of closed subspaces. This method was initially proposed by 
J. Von Neumann \cite{von1951functional} in the case of two closed subspaces and 
later extended by Halperin \cite{Halperin1962} in the case of $n$ closed subspaces. 
The procedure is summarized in Algorithm~\ref{alg:MAP}: starting from an arbitrary 
initial point $x_0$, it successively applies the projections $P_{M_1}, \ldots, P_{M_n}$ 
and repeats this sweep until the iterates stabilize. The convergence of this scheme 
is guaranteed by the following theorem, which states that the sequence of iterates 
converges to the orthogonal projection of $x_0$ onto the intersection 
$M = \bigcap_{i=1}^n M_i$.

\begin{algorithm}[H]
\caption{Method of Alternating Projections (Halperin, 1962)}
\label{alg:MAP}
\KwIn{Projections $P_{M_1}, \ldots, P_{M_n}$ onto closed subspaces of a Hilbert space $H$;
      initial point $x_0 \in H$;
      tolerance $\varepsilon > 0$;
      maximum number of iterations $K_{\max}$.}
\KwOut{Approximation $\hat{x} \approx P_M x_0$, where $M = \bigcap_{i=1}^n M_i$.}

$x \leftarrow x_0$\;
\For{$k \leftarrow 1$ \KwTo $K_{\max}$}{
    $y \leftarrow x$\;
    \For{$i \leftarrow 1$ \KwTo $n$}{
        $y \leftarrow P_{M_i}\, y$\;
    }
    \lIf{$\|y - x\| \leq \varepsilon$}{ \Return $y$\; }
    $x \leftarrow y$\;
}
\Return $x$\;
\end{algorithm}

In practice, the algorithm is frequently stopped when the displacement between consecutive 
iterates falls below a prescribed tolerance, i.e., when $\|y - x\| \leq \varepsilon$, but there may be other stopping criteria.

\begin{theorem}[Halperin, 1962]
    Let $M_1, M_2, \ldots, M_n$ be closed (linear) subspaces of a Hilbert space $H$ 
    with $M=\cap_{i=1}^n M_i$. Then,
    \[
    \lim_{k \to \infty} \|(P_{M_n}P_{M_{n-1}}\ldots P_{M_1})^k x - P_Mx\| = 0.
    \]
\end{theorem}

Before Halperin, Kaczmarz had proposed the MAP for the particular case of solving systems of linear equations \cite{karczmarz1937angenaherte}. There are variants of the Method of Alternating Projections for more general cases such as closed convex sets (e.g. \cite{gubin1967method, dykstra1983algorithm}). The MAP has a wide variety of applications; moreover, it has been discovered that many convergent iterative algorithms, widely used in numerical analysis, are in fact projection algorithms. For more information on the Method of Alternating Projections and its variants, see \cite{escalante2011alternating}. 

One characteristic of the MAP scheme is the slow convergence that can be observed in practice when the angles between the involved subspaces are small \cite{bauschke2003accelerating}. There have been many proposals to accelerate the method of alternating projections. For example, we can mention \cite{gubin1967method, gearhart1989acceleration, bauschke2003accelerating, censor2006computational}. More recently, Randomized Kaczmarz's ideas have attracted interest \cite{bai2018greedy, chen2023greedy, needell2010randomized, needell2014paved, needell2015randomized, strohmer2009randomized, zouzias2013randomized}. On the other hand, there have also recently been proposals for acceleration based on numerical optimization \cite{hernandez2011unconstrained, hernandez2005alternating}.

Iterations of the Method of Alternating Projections can be viewed as a fixed point scheme $x_k = Tx_{k-1}$, for the bounded linear non-expansive operator $T= P_{M_n}P_{M_{n-1}}\ldots P_{M_1}$ (i.e., $\|T\| \leq 1$, where $\|.\|$ is the linear-operator induced norm) \cite{bauschke2003accelerating, hernandez2011unconstrained, lopez2016acceleration, gearhart1989acceleration}. This scheme converges to the set of fixed points $Fix \; T = M = \cap_{i=1}^n M_i$, that is,
\[
\lim_{k \to \infty} \| T^kx - P_{Fix \; T}x\| = 0.
\]

In \cite{bauschke2003accelerating}, Bauschke et al. propose a generalization of the accelerated schemes in \cite{gubin1967method, gearhart1989acceleration}. This generalization is developed for non-expansive operators $T$ in general, including the operator $T=\frac{1}{n}\sum\limits_{i=1}^n P_{M_i}$ associated with the Cimmino's method \cite{cimmino1938}, another projection-based method generally slower than the method of alternating projections, but with better properties for implementation on parallel machines (see \cite{benzi2004gianfranco} for more information on the Cimmino's method).

In \cite{hernandez2011unconstrained}, a scheme accelerated by the Conjugate Gradient Method (CG) \cite{hestenes1952methods} is proposed for the case where the operator $T$ is self-adjoint. This scheme presents quite satisfactory performance in accelerating the Cimmino's method, but it cannot be applied to accelerate the Von Neumann-Halperin method unless it is done on a symmetrized, usually ill-conditioned, operator. 

In this article, for $H=\mathbb{R}^n$, we propose to use an acceleration scheme for the MAP operator based on the GMRES method \cite{saad1986gmres, saad2003iterative}. Although the GMRES method does not usually have better performance than CG, since it does not have to symmetrize the operator, better results are expected. We will also compare the proposed GMRES acceleration of the Method of Alternating Projections versus the Conjugate Gradients acceleration of the Cimmino method proposed in \cite{hernandez2011unconstrained} and the Randomized Kaczmarz accelerations.

\section*{Method of Alternating Projections Acceleration by Solving (Singular) Linear Systems}
In the case where the sets $M_i$ are subspaces of $\mathbb{R}^n$, each projection 
$P_{M_i}$ is a linear operator. Consequently, the operator
\[
T = P_{M_n}P_{M_{n-1}}\ldots P_{M_1}
\]
is also linear. Furthermore,
\[
Fix \; T = N(I-T) = M = \bigcap_{i=1}^n M_i,
\]
where $N(A)$ denotes the kernel of $A$. Feasible points $x \in M$ can be found by 
simply solving the homogeneous linear system
\[
Ax = 0, \qquad \text{where} \quad A = I - T.
\]
It is important to note that if there are non-trivial feasible solutions, then 
$A$ is a singular non-symmetric operator. The homogeneous system is always 
consistent, because there is always the trivial solution $x = 0$.

GMRES resolution of non-symmetric singular systems has been previously studied 
by Brown and Walker \cite{brown1997gmres}. This article establishes the conditions 
under which the method converges. The following theorem offers important results 
for our purposes: it guarantees that, when $N(A) = N(A^T)$, GMRES applied to 
$Ax = b$ does not break down before reaching a least-squares solution, and that 
the only possible breakdown afterwards is benign, occurring either through 
degeneracy of the Krylov subspace (when the system is consistent) or through 
rank deficiency of the small least-squares problem (when it is not). In our setting 
$b = 0$ and $A = I - T$, so the system is always consistent and the breakdown is 
necessarily of the first type. The algorithm that implements these ideas is 
summarized in Algorithm~\ref{alg:GMRES}.

\begin{theorem}[Brown et al.]
    Let $A \in \mathbb{R}^{n\times n}$. GMRES determines a least-squares solution of 
    $Ax=b$ without breakdown for all $b$ and $x_0$ if and only if $N(A)=N(A^T)$. 
    If $N(A)=N(A^T)$ and the least-squares solution is reached at some step, then GMRES 
    breaks down at the next step, with breakdown through degeneracy of the Krylov 
    subspace if $Ax=b$ is consistent and through rank deficiency of the least-squares 
    problem $\min_{z\in \mathcal{K}_k} \| b - A(x_0 + z)\|_2$ otherwise. Here, 
    $\mathcal{K}_k = \mathcal{K}_k(A, r_0) = \mathrm{span}\{r_0, A r_0, \ldots, 
    A^{k-1} r_0\}$ is the $k$th Krylov subspace determined by $A$ and the initial 
    residual $r_0 = b - A x_0$.
\end{theorem}

\begin{algorithm}[H]
\caption{GMRES for the singular system $(I-T)x=0$ (with modified Gram--Schmidt)}
\label{alg:GMRES}
\KwIn{Operator $T$ (accessible only through operator-vector products $Tv$);
      initial guess $x_0 \in \mathbb{R}^n$;
      tolerance $\varepsilon > 0$;
      maximum dimension of the Krylov subspace $m_{\max}$.}
\KwOut{Approximation $\hat{x} \approx P_{Fix \; T}\,x_0$.}

$r_0 \leftarrow T x_0 - x_0$\; \tcp*{initial residual: $r_0 = -(I-T)x_0$, since $b=0$}
$\beta \leftarrow \|r_0\|_2$\;
\lIf{$\beta \leq \varepsilon$}{\Return $x_0$\;}
$v_1 \leftarrow r_0 / \beta$\;
$H \leftarrow [\,]$\; \tcp*{upper Hessenberg matrix}
$V \leftarrow [v_1]$\;

\For{$j \leftarrow 1$ \KwTo $m_{\max}$}{
    $w \leftarrow v_j - Tv_j$\; \tcp*{$(I-T)v_j$; no matrix $A$ is built}
    \For{$i \leftarrow 1$ \KwTo $j$}{
        $h_{i,j} \leftarrow \langle w, v_i \rangle$\;
        $w \leftarrow w - h_{i,j}\, v_i$\; \tcp*{MGS: subtract one at a time}
    }
    $h_{j+1,j} \leftarrow \|w\|_2$\;
    \lIf{$h_{j+1,j} \leq \varepsilon$}{
        \textbf{break}\; \tcp*{benign breakdown: invariant subspace reached}
    }
    $v_{j+1} \leftarrow w / h_{j+1,j}$\;
    Append $v_{j+1}$ to $V$\;
    Solve $\min_{y \in \mathbb{R}^{j}} \|\beta e_1 - H_j\, y\|_2$ for $y_j$\; \tcp*{small least-squares problem}
    $x_j \leftarrow x_0 + V_j\, y_j$\;
    \lIf{$\|\beta e_1 - H_j\, y_j\|_2 \leq \varepsilon$}{
        \Return $x_j$\;
    }
}
\Return $x_j$\;
\end{algorithm}

In a finite-dimensional real vector space with an inner product, the adjoint 
operator of a linear operator $T$ is represented matrix-wise as the transpose 
matrix. In the Method of Alternating Projections on subspaces of $\mathbb{R}^n$, 
although $A = I - T$ is not necessarily a symmetric matrix, it is true that
\[
Fix \; T = N(A) = N(A^T) = Fix \; T^T,
\]
because the fixed points of
\[
T = P_{M_n}P_{M_{n-1}}\ldots P_{M_1}
\]
are the same fixed points of
\[
T^T = P_{M_1}P_{M_2}\ldots P_{M_n}.
\]
In both cases,
\[
Fix \; T = Fix \; T^T = M = \bigcap_{i=1}^n M_i
\]
(the intersection of sets is commutative). Thus, by applying GMRES to 
$(I-T)x = 0$, we should obtain an element of $M$ (a feasible solution). 
Now, let us prove that GMRES converges to the best approximation of $x_0$ in $M$.

GMRES iterations verify that
\[
x_m = x_0 + q_m(A)r_0,
\]
where $r_0 = b-Ax_0$ and $q_m(A)$ is a polynomial of degree $m$ in $A$. Since 
projections onto subspaces are linear operators, and $Fix \; T = M$ is a subspace, 
GMRES iterations to solve $(I-T)x = 0$ verify that
\[
P_{Fix \; T} x_m = P_{Fix \; T} x_0 + P_{Fix \; T} \; q_m(I-T)r_0.
\]
In a finite-dimensional space, if $T$ is a non-expansive operator, then
\[
Im(I-T)^T = N(T)^\perp = (Fix \; T)^\perp,
\]
but, for MAP we have
\[
(Fix \; T)^\perp = (Fix \; T^T)^\perp = Im(I-T).
\]
Consequently, since $q_m(I-T)r_0 \in Im(I-T)$ and $Im(I-T) \perp Fix \; T$, we have
\[
P_{Fix \; T} \, q_m(I-T)r_0 = 0.
\]
So, all GMRES iterations $x_m$ verify that
\[
P_{Fix \; T} x_m = P_{Fix \; T} x_0.
\]
As the method converges to $x_*$ belonging to $Fix \; T$, i.e. $P_{Fix \; T} x_* = x_*$. 
Hence, we conclude that
\[
x_* = P_{M} x_0.
\]
We can summarize the previous result in the following theorem:
\begin{theorem}\label{thm:main}
Let $T = P_{M_n}P_{M_{n-1}}\ldots P_{M_1}$, and let $M = \bigcap_{i=1}^n M_i$. 
Since $Fix \; T = Fix \; T^T = M$, we have $N(I-T) = N(I-T^T)$. Then GMRES 
applied to the system $(I-T)x = 0$ converges to $P_M x_0$, the best approximation 
of the initial iterate $x_0$ in $M$.
\end{theorem}
We must emphasize that to apply the GMRES method it is never necessary to construct 
the matrix $A=I-T$, but rather to apply the product 
$Ax = x - Tx = x - P_{M_n}P_{M_{n-1}}\ldots P_{M_1}x$ for some $x$ when necessary, 
by projecting consecutively on the subspaces involved. The resulting implementable 
version of the algorithm, in which $T$ is expanded in terms of the projections 
$P_{M_1}, \ldots, P_{M_n}$, is presented in Algorithm~\ref{alg:MAP-GMRES}.

\begin{algorithm}[H]
\caption{GMRES for the singular system $(I-T)x=0$, with 
$T = P_{M_n}P_{M_{n-1}}\ldots P_{M_1}$ (modified Gram--Schmidt)}
\label{alg:MAP-GMRES}
\KwIn{Projections $P_{M_1}, \ldots, P_{M_n}$ onto closed subspaces of $\mathbb{R}^n$;
      initial guess $x_0 \in \mathbb{R}^n$;
      tolerance $\varepsilon > 0$;
      maximum dimension of the Krylov subspace $m_{\max}$.}
\KwOut{Approximation $\hat{x} \approx P_{Fix \; T}\,x_0$.}

\tcp{initial residual: $r_0 = -(I-T)x_0 = T x_0 - x_0$}
$u \leftarrow x_0$\;
\For{$i \leftarrow 1$ \KwTo $n$}{
    $u \leftarrow P_{M_i}\, u$\;
}
$r_0 \leftarrow u - x_0$\;
$\beta \leftarrow \|r_0\|_2$\;
\lIf{$\beta \leq \varepsilon$}{\Return $x_0$\;}
$v_1 \leftarrow r_0 / \beta$\;
$H \leftarrow [\,]$\; \tcp*{upper Hessenberg matrix}
$V \leftarrow [v_1]$\;

\For{$j \leftarrow 1$ \KwTo $m_{\max}$}{
    \tcp{compute $(I-T)v_j = v_j - T v_j$ by applying the projections}
    $u \leftarrow v_j$\;
    \For{$i \leftarrow 1$ \KwTo $n$}{
        $u \leftarrow P_{M_i}\, u$\;
    }
    $w \leftarrow v_j - u$\; \tcp*{$w = (I-T)v_j$; no matrix $A$ is built}
    \For{$i \leftarrow 1$ \KwTo $j$}{
        $h_{i,j} \leftarrow \langle w, v_i \rangle$\;
        $w \leftarrow w - h_{i,j}\, v_i$\; \tcp*{MGS: subtract one at a time}
    }
    $h_{j+1,j} \leftarrow \|w\|_2$\;
    \lIf{$h_{j+1,j} \leq \varepsilon$}{
        \textbf{break}\; \tcp*{benign breakdown: invariant subspace reached}
    }
    $v_{j+1} \leftarrow w / h_{j+1,j}$\;
    Append $v_{j+1}$ to $V$\;
    Solve $\min_{y \in \mathbb{R}^{j}} \|\beta e_1 - H_j\, y\|_2$ for $y_j$\; \tcp*{small least-squares problem}
    $x_j \leftarrow x_0 + V_j\, y_j$\;
    \lIf{$\|\beta e_1 - H_j\, y_j\|_2 \leq \varepsilon$}{
        \Return $x_j$\;
    }
}
\Return $x_j$\;
\end{algorithm}

\section*{Numerical Experiments}
We will compare the proposed acceleration scheme with some well-known acceleration 
schemes that can be found in the literature. Numerical experiments will be performed 
with $A$-projections (projections where the orthogonality condition is verified for 
the inner product $\langle.,.\rangle_A$), which is equivalent to solving a saddle 
point problem \cite{hernandez2005alternating, hernandez2011unconstrained}. Linear systems of saddle point type 
arise in a wide variety of applications throughout computational science and 
engineering \cite{benzi2005numerical}. 

\subsection*{Saddle point problems resolution by MAP}
We are interested in solving the system of linear equations:
\begin{equation}\label{saddle}
    \begin{pmatrix}
A & B^T \\
B & 0
\end{pmatrix}
\begin{pmatrix}
    x \\
    \lambda
\end{pmatrix}
=
\begin{pmatrix}
    f \\
    0
\end{pmatrix},
\end{equation}
where
\begin{description}
    \item[(H1)] $B \in \mathbb{R}^{m \times n}$ is a matrix with $\operatorname{rank}(B)= m$, $m \leq n$,
    \item[(H2)] $A \in \mathbb{R}^{n \times n}$ is a symmetric positive definite matrix.
\end{description}
The next proposition is well-known \cite{benzi2005numerical, hernandez2005alternating}. Nevertheless we present it for the sake of completeness.
\begin{propo}\label{BAP}
    Under hypotheses $(H1, H2)$, if $(x, \lambda)^T$ solves \eqref{saddle}, then the 
    vector $x$ corresponds to the orthogonal projection of $x_u = A^{-1}f$ onto 
    $N(B)$, in the inner product $\langle .,.\rangle_A$.
\end{propo}

Let us partition the matrix $B\in \mathbb{R}^{m \times n}$ into $r$ row blocks:
\[
B^T = [B_1^T, B_2^T, \ldots, B_r^T].
\]
Then,
\[
N(B) = \bigcap_{i=1}^r N(B_i).
\]
From Proposition \ref{BAP}, it follows that the solution $x$ for the saddle point 
problem \eqref{saddle} is the $A$-orthogonal projection of $x_u=A^{-1}f$ onto the 
subspace $N(B) = \bigcap_{i=1}^r N(B_i)$. Hence, this problem can be viewed as a 
Best Approximation Problem (BAP) in the inner product $\langle.,.\rangle_A$ and can 
be solved by the Method of Alternating Projections \cite{hernandez2011unconstrained}. 
As each $N(B_i)$ is also a subspace, we can obtain the $A$-projection of any vector 
$y \in \mathbb{R}^n$ onto $N(B_i)$, $P_{M_i}y$, by computing $x$ from the smaller 
saddle point problem:
\begin{equation*}\label{saddle2}
    \begin{pmatrix}
A & B_i^T \\
B_i & 0
\end{pmatrix}
\begin{pmatrix}
    x \\
    \lambda
\end{pmatrix}
=
\begin{pmatrix}
    Ay \\
    0
\end{pmatrix}.
\end{equation*}
These smaller saddle point problems can be solved by any method, whose suitability 
depends on the size of each block.

\subsection*{Results}
Different accelerated versions of projection methods will be evaluated in solving 
saddle point problems from a selection of problems in the CUTEr collection 
\cite{gould2001general}. 

Numerical experiments were performed on Octave 9.1.0 running on a computer with an 
AMD Ryzen 7 5700G processor. In all the experiments and for all the considered 
methods, we stop the process when the norm-$A$ of the residual is less than or 
equal to $10^{-6}$. 

For all sets of problems, the matrix $A \in \mathbb{R}^{n\times n}$ is a symmetric 
and positive definite matrix, and $B \in \mathbb{R}^{m\times n}$ is a full row 
rank matrix. It is important to note that an augmented Lagrangian technique has 
been used to improve the behavior of all considered methods, and also to guarantee 
that the $A$ matrix is positive definite. To be precise, the $A$ matrix in 
\eqref{saddle} has been replaced by
\[
A + \gamma\, B^T B, \qquad \text{where} \quad \gamma = \frac{\|A\|_F}{\|B\|^2_F}
\]
\cite{benzi2005numerical}. Moreover, we also added to the left upper block of 
\eqref{saddle} the scalar $0.01$ times the identity matrix, to improve the 
condition number of that block.

The source code for the experiments can be found in the following GitHub repository: 
\href{https://github.com/GorkaEP/Saddle-Point-Problem}{https://github.com/GorkaEP/Saddle-Point-Problem}.

Numerical comparisons with the proposed acceleration scheme were divided into two groups:
\begin{itemize}
    \item \textbf{Group 1.} Comparisons with schemes based on the von Neumann-Halperin 
    operator or its symmetrization. Here we will compare the proposed acceleration with:
    \begin{itemize}[parsep=0.1cm]
        \item Classical von Neumann-Halperin method (MAP).
        \item Conjugate gradient-based acceleration to minimize 
        $f(x) = \frac{1}{2} \|x-Tx\|^2$ proposed in \cite{hernandez2011unconstrained}.      
    \end{itemize}
    \item \textbf{Group 2.} Comparisons of the proposed acceleration with:
    \begin{itemize}
        \item Conjugate gradient-based acceleration to solve the singular symmetric 
        problem $(I - T_S)x = 0$, where $T_S=\frac{1}{n}\sum\limits_{i=1}^n P_{M_i}$ 
        is the operator associated with the Cimmino's method \cite{hernandez2011unconstrained}.
        \item The Greedy Randomized Kaczmarz method proposed in \cite{bai2018greedy}.
    \end{itemize}
\end{itemize}

\subsubsection*{Results for the Group 1}

\begin{figure} 
        \centering    
        \includegraphics[width=0.8\textwidth]{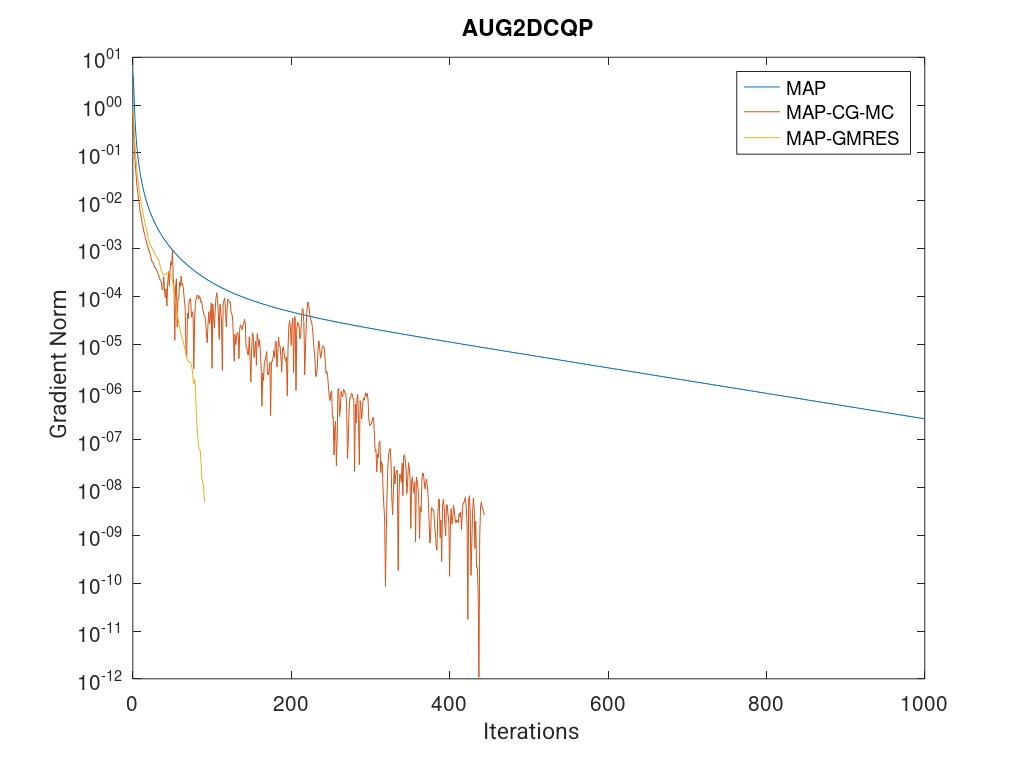}
        \caption{Comparison of group 1 methods on AUG2DCQP (CUTEr).}
        \label{fig:AUG2DCQPgraph}
    \end{figure}

Table \ref{tab:vonNeumannHalperinMethods} presents the results of the comparison of 
the proposed GMRES scheme with methods based on the von Neumann-Halperin operator 
(or its symmetrization) applied to a selection of saddle-point problems from the 
CUTEr library. As can be seen, the acceleration proposal based on GMRES has a higher 
performance than the classical MAP (MAP) and its acceleration applying conjugate 
gradients to the symmetrized operator \cite{hernandez2011unconstrained} (MAP-CG-MC) 
in almost all the problems solved, both in number of iterations and in execution 
time. This superior performance can also be observed for the AUG2DCQP problem in 
Figure \ref{fig:AUG2DCQPgraph}.

    \begin{table}
        \centering
        \begin{tabular}{|l|l|l|l|l|l|l|l|l|}
            \midrule
            \multicolumn{3}{|c|}{Problem} &  \multicolumn{2}{|c|}{MAP} & \multicolumn{2}{|c|}{MAP-CG-MC}  & \multicolumn{2}{|c|}{MAP-GMRES}\\ \hline
            Name & m & n & Iter & Time & Iter & Time & Iter & Time \\ \hline
            PRIMALC2 & 7 & 238 & 6 & 0.029283 & 4 & 0.041855 & \textbf{4} & \textbf{0.025755} \\ \hline
            PRIMALC1 & 9 & 239 & 6 & 0.03845 & 5 & 0.062401 & \textbf{5} & \textbf{0.040444} \\ \hline
            PRIMAL4 & 75 & 1564 & 20 & 106.7301 & 14 & 181.2912 & \textbf{11} & \textbf{78.79399} \\ \hline
            PRIMAL1 & 85 & 410 & 25 & 7.060531 & 14 & 7.876206 & \textbf{13} & \textbf{4.084433} \\ \hline
            PRIMAL2 & 96 & 745 & 23 & 20.18988 & 14 & 23.66418 & \textbf{12} & \textbf{11.04813} \\ \hline
            PRIMAL3 & 111 & 856 & 41 & 52.91835 & 19 & 54.17642 & \textbf{16} & \textbf{22.17518} \\ \hline
            DUALC1 & 215 & 223 & 3 & 0.139782 & 3 & 0.283966 & \textbf{1} & \textbf{0.097391} \\ \hline
            DUALC5 & 278 & 285 & 4 & 0.279207 & 4 & 0.564005 & \textbf{2} & \textbf{0.208993} \\ \hline
            DUALC8 & 503 & 510 & 2 & \textbf{0.45704} & 3 & 1.380096 & \textbf{1} & 0.457258 \\ \hline
            QGROW15 & 300 & 645 & 8 & 1.662005 & 6 & 2.568962 & \textbf{5} & \textbf{1.242921} \\ \hline
            QGROW22 & 440 & 946 & 8 & 3.617278 & 9 & 8.244714 & \textbf{6} & \textbf{3.316471} \\ \hline
            MOSARQP1 & 700 & 3200 & 3 & 5.786842 & 4 & 15.71982 & \textbf{2} & \textbf{5.777938} \\ \hline
            QSCFXM3 & 990 & 1800 & - & - & 293 & 3332.4 & \textbf{102} & \textbf{480.5315} \\ \hline
            KSIP & 1001 & 1021 & 296 & 530.1636 & 21 & 71.82051 & \textbf{31} & \textbf{54.98155} \\ \hline
            AUG2DCQP & 1600 & 3280 & - & - & 445 & 4407.819 & \textbf{92} & \textbf{460.656} \\ \hline
            CVXQP1\_S & 50 & 100 & 671 & 4.1122 & 54 & 0.666828 & \textbf{32} & \textbf{0.215098} \\ \hline
            CVXQP2\_S & 25 & 100 & 23 & 0.06625 & 18 & 0.100785 & \textbf{11} & \textbf{0.035338} \\ \hline
            CVXQP3\_S & 75 & 100 & - & - & 200 & 4.095119 & \textbf{68} & \textbf{0.740479} \\ \hline
            \multicolumn{9}{p{\linewidth}}{\footnotesize Note: the symbol '-' indicates that the method did not converge at the maximum number of iterations (1000).} \\
        \end{tabular}
        \caption{Results of the group 1 of experiments.}
        \label{tab:vonNeumannHalperinMethods}
    \end{table}
    
    \subsubsection*{Results for the Group 2}

\begin{figure}
        \centering    
        \includegraphics[width=0.8\textwidth]{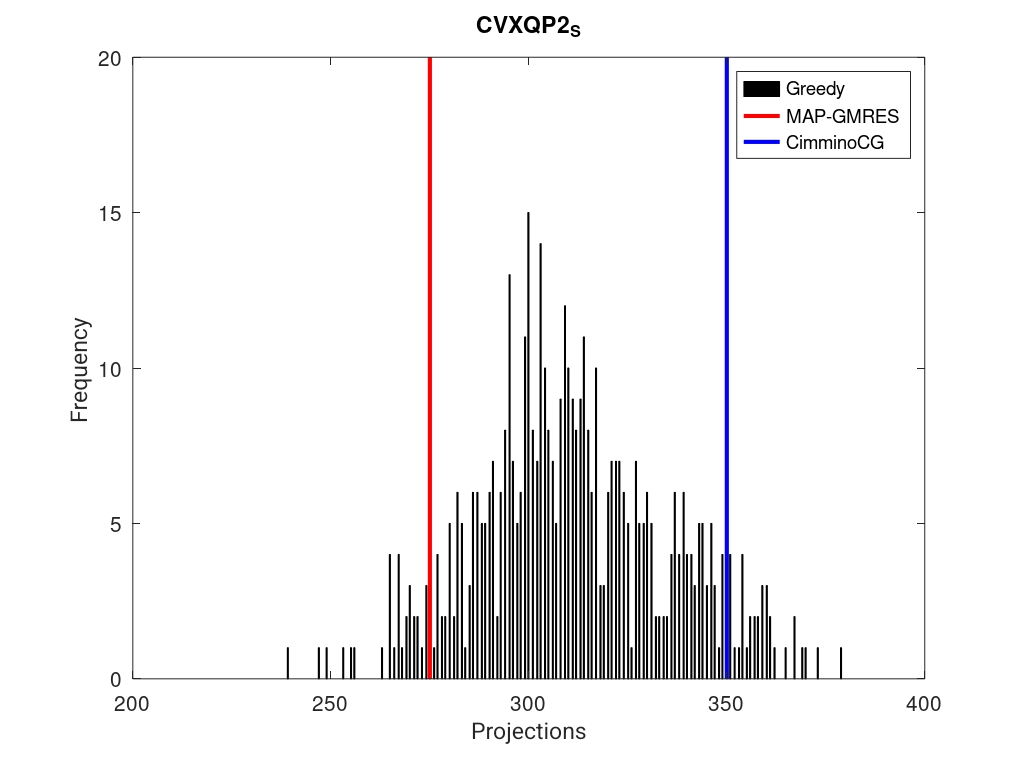}
        \caption{Projections needed for convergence (CVXQP2\_S).}
        \label{fig:histogram}
    \end{figure}

\begin{figure}
        \centering    
        \includegraphics[width=0.8\textwidth]{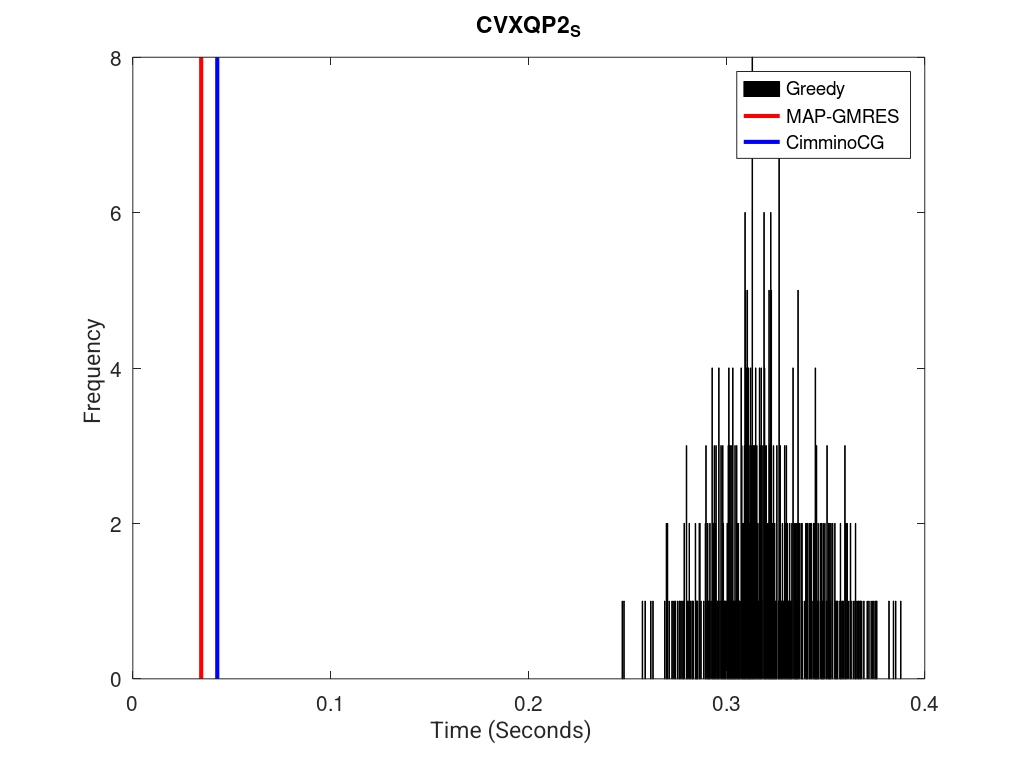}
        \caption{CPU time required for convergence (CVXQP2\_S).}
        \label{fig:histogramtime}
    \end{figure}

Table \ref{tab:MAPMethods} compares the proposed acceleration of MAP by GMRES 
(MAP-GMRES) with the conjugate gradient acceleration of the Cimmino method 
(Cimmino-GC) proposed in \cite{hernandez2011unconstrained} and with the Kaczmarz 
Greedy method proposed in \cite{bai2018greedy}. The table details the total number 
of projections used for each method and the time taken. In the case of the Greedy 
Kaczmarz method, since these values are random, averages over $500$ runs are shown. 
While there is no clear winner in terms of the number of projections used, in terms 
of average runtime, the Cimmino-GC and MAP-GMRES methods performed best, with no 
clear winner between the latter two. As an example, Figure \ref{fig:histogram} shows, 
for a problem from the CUTEr library (CVXQP2\_S), the number of projections required 
for convergence, while Figure \ref{fig:histogramtime} shows the execution time for 
the same problem ($500$ runs).

    \begin{table}
        \centering
        \begin{tabular}{|l|l|l|l|l|l|l|l|l|}
            \midrule
            \multicolumn{3}{|c|}{Problem} &    \multicolumn{2}{|c|}{ Cimmino-CG}&\multicolumn{2}{|c|}{MAP-GMRES} &  \multicolumn{2}{|c|}{Kaczmarz Greedy} \\ \hline
            Name & m & n &   Proj&Time &Proj& Time & Proj*& Time  \\ \hline
            PRIMALC2& 7 & 238 &   \textbf{21}&0.03145&28 & \textbf{0.02575}& 40.564 & 0.05157\\ \hline
            PRIMALC1& 9 & 239 &   \textbf{27}&\textbf{0.02839}&45 & 0.04044& 54.792 & 0.08188\\ \hline
            PRIMAL4& 75 & 1564 &   \textbf{750}&\textbf{57.6934}&825 & 78.7939& 738.506 & 59.6576\\ \hline
            PRIMAL1& 85 & 410 &   850&\textbf{2.76655}&1105 & 4.08443& \textbf{684.014} & 4.98433\\ \hline
            PRIMAL2& 96 & 745 &   960&\textbf{8.31414}&1152 & 11.0481& \textbf{766.85} & 11.6874\\ \hline
            PRIMAL3& 111 & 856 &   1443&\textbf{18.7870}&1776 & 22.1751& \textbf{1404.71} & 28.3562\\ \hline
            DUALC1& 215 & 223 &   645&0.13419&\textbf{215} & \textbf{0.09739}& 217 & 1.29577\\ \hline
            DUALC5& 278 & 285 &   1112&0.26150&556 & \textbf{0.20899}& \textbf{247.152} & 1.93922\\ \hline
            DUALC8& 503 & 510 &   1509&0.62405&\textbf{503} & \textbf{0.45725}& 509.098 & 8.40810\\ \hline
            QGROW15& 300 & 645 &   1500&\textbf{1.10462}&1500 & 1.24292& \textbf{1377.862} & 15.5873\\ \hline
            QGROW22& 440 & 946 &   3080&3.58143&2640 & \textbf{3.31647}& \textbf{2104.392} & 41.5511\\ \hline
            MOSARQP1& 700 & 3200 &   1400&\textbf{4.84117}&1400 & 5.77793& \textbf{1385.738} & 87.5882\\ \hline 
            QSCFXM3& 990& 1800& \textbf{90090}& \textbf{430.408}& 100980& 480.531& -&-\\ \hline 
            KSIP& 1001& 1021& \textbf{6006}& \textbf{12.0417}& 31031& 54.9815& -&-\\ \hline 
            AUG2DCQP& 1600& 3280& \textbf{102400}& \textbf{322.168}& 147200& 460.656& -&-\\ \hline
            CVXQP1\_S& 50 & 100 &   2550&0.30918&\textbf{1600} & \textbf{0.21509}& 12164.84 & 20.0816\\ \hline
            CVXQP2\_S& 25 & 100 &   350&0.04264&\textbf{275} & \textbf{0.03533}& 311.56 & 0.31929\\ \hline
            CVXQP3\_S& 75& 100& 8925& 1.24346& \textbf{5100}& \textbf{0.74047}& -&-\\ \hline
            \multicolumn{9}{p{\linewidth}}{\footnotesize * Mean of number of projections.} \\
        \end{tabular}
        \caption{Results of the group 2 of experiments.}
        \label{tab:MAPMethods}
    \end{table}

\section*{Final remarks}
Our preliminary results show that the GMRES acceleration of the Method of Alternating 
Projections outperforms the idea of symmetrizing the operator to apply Conjugate 
Gradients proposed in \cite{hernandez2005alternating,hernandez2011unconstrained}. The experiments also show 
better runtime performance than the Kaczmarz Greedy algorithm \cite{bai2018greedy} 
in most cases. Regarding the Cimmino acceleration using conjugate gradients, the 
MAP-GMRES proposal was competitive for some problems, although it must be taken 
into account that the Cimmino algorithm is by its nature more parallelizable than 
the method of alternating projections.

\section*{Acknowledgments}
The authors would like to thank Marcos Raydan for interesting discussions and 
constructive suggestions on this topic and for a careful reading of this document.

\bibliographystyle{plain}
\bibliography{Biblio}

\end{document}